\documentclass{elsart}

\usepackage{natbib}
\usepackage{amssymb}

\newtheorem{theorem}{Theorem}
\newtheorem{lemma}{Lemma}
\newtheorem{corollary}{Corollary}
\begin{document}
\begin{frontmatter}


 \title{On posterior distribution of Bayesian wavelet thresholding}%\thanksref{label1}}
 %\thanks[label1]{111}
 \author{Heng Lian}
 \ead{henglian@ntu.edu.sg}
  \address{Division of Mathematical Sciences, SPMS,\\ Nanyang Technological University, \\Singapore 637616}

%\title{}


%\author{}

%\address{}

\begin{abstract}
We investigate the posterior rate of convergence for wavelet
shrinkage using a Bayesian approach in general Besov spaces.
Instead of studying the Bayesian estimator related to a particular
loss function, we focus on the posterior distribution itself from
a nonparametric Bayesian asymptotics point of view and study its
rate of convergence. We obtain the same rate as in
\citet{abramovich04} where the authors studied the convergence of
several Bayesian estimators.

\end{abstract}


\end{frontmatter}

% main text
\section{Introduction}\label{sec:introduction}

Infinite-dimensional Bayesian methods have become quite popular
recently, due to both the computational and theoretical advances
in this field. There are many results concerning posterior
convergence using appropriate priors. These developments originate
from the consideration of density estimation problems. In these
problems, given the prior $\Pi_n$ on the set $\mathcal{P}$ of
probability distributions, the posterior is a random measure:
\begin{eqnarray}\label{eqn:posterior}
\Pi_n(B|X_1,\ldots,X_n)=\frac{\int_B\prod_{i=1}^np(X_i)d\Pi_n(P)}{\int\prod_{i=1}^np(X_i)d\Pi_n(P)}
\end{eqnarray}
We say that the posterior is consistent if
\begin{eqnarray*}\label{eqn:rate}
\Pi_n(P\in \mathcal{P}: d(P,P_0)>\epsilon|X)\rightarrow 0 \mbox{
in $P_0^n$ probability}.
\end{eqnarray*}
where $P_0$ is the true distribution and $d$ is some suitable
distance function between probability measures.

To study rates of convergence, let $\epsilon_n$ be a sequence
decreasing to zero, we say the rate is at least $\epsilon_n$ if
for sufficiently large constant $M$
\[
\Pi_n(P: d(P,P_0)\ge M\epsilon_n|X)\rightarrow 0 \mbox{ in $P_0^n$
probability}.
\]

It turns out the convergence rates are closely related to the
existence of tests that separate the hypotheses in convex sets.

The most general result appeared recently in \cite{ghosal07},
where the formulation includes both density estimation and
regression problems, and the results also extend to non-iid cases
such as stationary and nonstationary sequence of observations. In
this general context, the definition of the posterior convergence
rate is similar except the measure $P$ and $P_0$ represent the
distribution on the data and thus depends on sample size $n$, and
the observations are no longer i.i.d. so the likelihood used in
(\ref{eqn:posterior}) must be changed to a more general form
accordingly.

Another relatively recent development in statistics is the
investigation of wavelet method which has found numerous
applications in engineering as well. There are many theoretical
results explaining why wavelet transformation is effective, from
both the frequentist and the Bayesian point of view. These
well-known results include the now widely celebrated works of
David Donoho and his collaborators \citep{donoho94,donoho96}. The
property that distinguishes these works from previous results is
that a single estimator can achieve the minimax rate over a range
of function spaces including functions with inhomogeneous
smoothness, whose minimax rate cannot be achieved by the simpler
linear estimator. The sparsity of the coefficients for the
function in an appropriate basis is the key to the success of the
wavelet thresholding approach.

Bayesian approach to function estimation in Besov spaces has been
investigated in \cite{abramovich98,abramovich04}. In these
approaches, after specifying an appropriate prior, the Bayesian
estimator is obtained from the posterior and investigated from the
frequentist point of view. In particular, they study the rate of
convergence of different point estimators including the posterior
mean and posterior median as well as other estimators derived from
the posterior distribution. The theoretical results in
\cite{abramovich04} show that some Bayesian estimators can achieve
the better-than-linear rates if an appropriate prior is chosen
that implicitly implements shrinkage or thresholding rule similar
to the frequentist approach.

In a Bayesian framework, most researchers are more interested in
the posterior as a distribution, instead of the point estimates
derived from specific loss function. The convergence of the
posterior distribution in this context has not been studied. This
paper intends to fill this gap. Using the same prior as in
\cite{abramovich98,abramovich04}, we show that the posterior
distribution has the same convergence rate as the point estimators
proposed in those papers.

We describe the model and present the main theorem in Section
\ref{sec:mainresult}. Some possible extensions for our result are
discussed in the final section.

\section{Main result}\label{sec:mainresult}
Consider the white noise model
\begin{eqnarray}\label{eqn:model}
dX(t)=f(t)dt+\sigma_n dW(t)
\end{eqnarray}
where $\sigma_n^2=1/n, f\in B_{p,q}^s[0,1]$ and $W$ is the
standard Brownian motion. Using wavelet basis on [0,1] with
sufficient regularity, the function $f$ can be expanded as
\begin{eqnarray*}
f&=&\sum_{k=0}^{2^{j_0}-1}\alpha_{j_0k}\phi_{j_0k}+\sum_{j\ge
j_0}\sum_{k=0}^{2^j-1}\beta_{jk}\psi_{jk}
\end{eqnarray*}
where $\phi_{j_0k}$ are the scaling functions and $\psi_{jk}$ are
the mother wavelets at resolution $j$, and $j_0$ is the lowest
resolution in the expansion. We assume $j_0=0$ for simplicity of
notation below.

The Besov spaces include the well-known Sobolev and H\"older
classes of function and also nearly contains the space of
functions of bounded variation. The norm for the Besov space with
parameter $s>\max(0,1/p-1/2), 1\le p\le \infty,$ and $1\le q\le
\infty$ is defined as
\begin{eqnarray*}
||f||_{B^s_{p,q}}&=&||P_0(f)||_{L^p}+(\sum_{j\ge
0}(2^{js}||Q_j(f)||_{L^p})^q)^{1/q}
\end{eqnarray*}
where $P_0(f)=\alpha_{00}\phi_{00}$ is the projection of $f$ on
the ``approximation space", and
$Q_j(f)=\sum_{k=0}^{2^j-1}\beta_{jk}\psi_{jk}$ is the projection
of $f$ onto the ``detail space".

In terms of the coefficients in the wavelet expansion, the Besov
norm can be equivalently defined by
\begin{eqnarray*}
||f||_{B^s_{p,q}}&\asymp
&||\beta||_{B^s_{p,q}}=|\alpha_{00}|+\left\{\sum_{j=0}^\infty
2^{j(s+1/2-1/p)q}||\beta_{j.}||_p^q\right\}^{1/q}
\end{eqnarray*}
Note that for cases where $q=\infty$ the usual change to the sup
norm is required. By abuse of notation, we also define $P_J\beta$
to be the sequence $\beta'$ such that $\beta'_{jk}=\beta_{jk}$
when $j\le J$ and $\beta'_{jk}=0$ when $j>J$.

The white noise model (\ref{eqn:model}) is closely related to the
nonparametric regression model \citep{brown96,donoho95b}:
\begin{eqnarray*}
Y_i&=&f\left(\frac{i}{n}\right)+z_i
\end{eqnarray*}
with standard normal noise. We choose to work with
(\ref{eqn:model}) for its simplicity of formulation.

After wavelet transformation for (\ref{eqn:model}), we get the
Gaussian sequence model:
\begin{eqnarray*}
X_{00}&=&\alpha_{00}^0+z_{00}/\sqrt{n}\\
X_{jk}&=&\beta_{jk}^0+z_{jk}/\sqrt{n}, j\ge 0, k=0,1,\ldots, 2^j-1
\end{eqnarray*}
where the superscript $0$ indicates the true parameter.

Using Bayesian approach for Gaussian sequence estimation, we put a
prior on $\beta_{jk}^0$:
\begin{eqnarray}\label{prior}
\beta_{jk}^0&\sim& \pi_jN(0,a_j^2)+(1-\pi_j)\delta_0
\end{eqnarray}
with hyperparameters $a_j^2\asymp 2^{-\alpha j}, \pi_j\asymp
2^{-\gamma j}$, for some $\alpha>1$, and $\gamma\ge 0$. This prior
is proposed in \citet{abramovich98}, and \citet{abramovich04}
investigated the optimality of some Bayesian estimators with this
prior. The choice of the hyperparameters must satisfy some
conditions for the prior to put positive mass on $B^s_{p,q}$
\citep[Theorem 1]{abramovich98}, although this is not our focus
here. In the following we assume $\alpha$ and $\gamma$ satisfy
these conditions. We also assume the value $\alpha_{00}^0$ is
known for simplicity, which does not affect our asymptotic result.

We also consider the sieve prior as in \citet{shenwong01}, and
define the prior $\Pi_n$ by $\Pi_n(A)=\sum_m \lambda_m\Pi_n^m(A)$
where $\lambda_m\propto 2^{-\mu m}$ for some $\mu>0$, and
$\Pi_n^m$ is a prior on $\beta_{jk}$ such that $\beta_{jk}\sim
N(0,2^{-\alpha j})$ when $j\le m$ and $\beta_{jk}=0$ when $j>m$.

The main result we obtain in this paper is the following:
\begin{theorem}\label{maintheorem}
Consider a bounded subset of the Besov space:
$B^s_{p,q}(B)=\{\beta\in B^s_{p,q}[0,1],
||\beta||_{B^s_{p,q}}<B\}$ and $\beta^0\in B^s_{p,q}(B)$. Let
$\alpha=2s+1$ for $p\ge 2$, and $\alpha=(2s+2-2/p)$ for $1\le
p<2$. Then for sufficiently large constant $M$, under the prior
(\ref{prior}), we have
\begin{eqnarray*}
\Pi_n\left(\{\beta_{jk}\}:
\sum_{j,k}(\beta_{jk}-\beta^0_{jk})^2>M\epsilon_n^2|X_{jk}\right)\rightarrow
0 \mbox{ in probability},
\end{eqnarray*}
where $\epsilon_n^2=(\log n)^2 n^{-2s/(2s+1)}$ when $p\ge 2$, and
$\epsilon_n^2=(\log n)^2 n^{-(2s+1-2/p)/(2s+2-2/p)}$ when $1\le
p<2$.
\end{theorem}

Remark: The above rate of convergence is the same as in
\citet{abramovich04} for posterior mean and posterior median,
except an extra $\log$ factor in our case, which we think might be
an artifact of our proofs.

\textbf{Proof of Theorem 1.} In the proof, we use $C$ to denote
generic constant whose value can change in difference locations.
We make use of the general result for Bayesian posterior rate of
convergence (Theorem 6 in \citet{ghosal07}), although we only use
a simpler version which corresponds to Theorem 2.1 in
\citet{ghosal00} in the iid case. Two conditions for the theorems
must be verified:

(I) $\log D(\epsilon_n, B^s_{p,q}(B),||.||_2)\le n\epsilon_n^2,$
where $D(\epsilon, F, ||.||)$ is the $\epsilon$-covering number of
the space $F$ with norm $||.||$.

(II) $\Pi_n^B(\beta\in B^s_{p,q}(B): ||\beta-\beta^0||_2^2\le
\epsilon_n^2)\ge \exp \{-Cn\epsilon_n^2\}$, where $\Pi_n^B$ denote
the prior distribution as in (\ref{prior}) constrained on
$B^s_{p,q}(B)$ by renormalization.

Corollary 2 in \citet{nickl07} gives the bracketing entropy number
for Besov spaces as $H_{B}(\epsilon,B^s_{p,q}(B),||.||_2)\precsim
\epsilon^{-1/s}$. Since bracketing entropy number is an upper
bound for usual entropy, $\epsilon_n$ defined in the statement of
the theorem obviously satisfies condition (I).

Condition (II) is verified as follows:

Since $\beta^0\in B^s_{p,q}(B)$, there exists $\delta$ such that
$||\beta^0||_{B^s_{p,q}}\le B-\delta$. Let $J=(\log_2 n)/\alpha$.
We have
\begin{eqnarray*}
&&\Pi_n^B(||\beta-\beta^0||_2^2\le\epsilon_n^2) \\
&\ge&\Pi_n(||\beta-\beta^0||_2^2\le\epsilon_n^2,
||\beta||_{B^s_{p,q}}<B\}\\
&\ge&\Pi_n(\sum_{j=0}^J\sum_k(\beta_{jk}-\beta_{jk}^0)^2\le\epsilon^2_n/2,
||P_J\beta||_{B^s_{p,q}}<B-\delta/2)\cdot\\
&&\Pi_n(\sum_{j=J+1}^\infty\sum_k(\beta_{jk}-\beta_{jk}^0)^2\le\epsilon^2_n/2,
||\beta-P_J\beta||_{B^s_{p,q}}<\delta/2)
\end{eqnarray*}
The above two terms are dealt with in the following two lemmas,
which provide a lower bound of $e^{-Cn\epsilon_n^2}$ and the
theorem is proved. $\quad\Box$

\begin{lemma}\label{lemma:secterm}
$\Pi_n(\sum_{j=J+1}^\infty\sum_k(\beta_{jk}-\beta_{jk}^0)^2\le\epsilon^2_n/2,
||\beta-P_J\beta||_{B^s_{p,q}}<\delta/2)$ is bounded away from
$0$.
\end{lemma}
\textbf{Proof.} Since $\sum_{j>J}\sum_k(\beta^0_{jk})^2\le
\sum_{j>J}C2^{-2js'}\le \epsilon_n^2/8$, where $s'=s$ for $p\ge 2$
and $s'=s+1/2-1/p$ otherwise, we have
\begin{eqnarray*}
&&\Pi_n(\sum_{j=J+1}^\infty\sum_k(\beta_{jk}-\beta_{jk}^0)^2\le\epsilon^2_n/2)\\
&\ge&\Pi_n(\sum_{j>J,k}\beta_{jk}^2\le\epsilon^2_n/8)\\
&\ge&1-8E[\sum_{j>J,k}\beta_{jk}^2]/\epsilon_n^2\\
&\ge&1-C\cdot 2^{-(\alpha-1)j}/\epsilon_n^2\\
&\rightarrow&1
\end{eqnarray*}

On the other hand,
$\Pi_n(||\beta-P_J\beta||_{B^s_{p,q}}<\delta/2)\ge\Pi_n(||\beta||_{B^s_{p,q}}<\delta/2)=:t>0$
when $\alpha$ and $\beta$ are chose appropriately such that
$\Pi_n(B^s_{p,q})>0$ (this is possible by \citet{abramovich98}).

Thus
$\Pi_n(\sum_{j=J+1}^\infty\sum_k(\beta_{jk}-\beta_{jk}^0)^2\le\epsilon^2_n/2,
||\beta-P_J\beta||_{B^s_{p,q}}<\delta/2)\rightarrow t>0$ as
$n\rightarrow \infty$. $\quad\Box$

\begin{lemma}\label{lemma:firstterm}
\[
\Pi_n(\sum_{j=0}^J\sum_k(\beta_{jk}-\beta_{jk}^0)^2\le\epsilon^2_n/2,
||P_J\beta||_{B^s_{p,q}}<B-\delta/2)\ge e^{-Cn\epsilon_n^2}
\]
\end{lemma}

\textbf{Proof.} This probability can be bounded from below using
the techniques in Section 5 of \citet{shenwong01}.

First we show
\begin{eqnarray}\label{temp}
&&\Pi_n(\sum_{j=0}^J\sum_k(\beta_{jk}-\beta_{jk}^0)^2\le\epsilon^2_n/2,
||P_J\beta||_{B^s_{p,q}}<B-\delta/2)\nonumber\\
&\ge&\Pi_n(\sum_{j=0}^J\sum_k(\beta_{jk}-\beta_{jk}^0)^2\le
c^2 \tau^2_n/\log n)
\end{eqnarray}
for a small enough constant $c$, where
$\tau_n=n^{-(s+1/2-1/p)/(2s+1)}$ when $p\ge 2$ and
$\tau_n=n^{-s/(2s+2-2/p)}$ when $1\le p <2$. Notice we obviously
have $\tau_n=o(\epsilon_n)$.

Case 1: $p\ge 2$.

Note $||\beta_{j.}||_p\le ||\beta_{j.}||_2$ when $p\ge 2$.
Conditioned on the event
$\sum_{j=0}^J\sum_k(\beta_{jk}-\beta_{jk}^0)^2\le c^2
\tau^2_n/\log n$, the $B^s_{p,q}$ norm for $P_J\beta-P_J\beta^0$
can be bounded as follows:
\begin{eqnarray*}
&&||P_J\beta-P_J\beta^0||_{B^s_{p,q}}\\
&=&(\sum_{j=0}^J 2^{j(s+1/2-1/p)q}||\beta_{j.}-\beta_{j.}^0||_p^q)^{1/q}\\
&\le&(\sum_{j=0}^J 2^{j(s+1/2-1/p)q}||\beta_{j.}-\beta_{j.}^0||_2^q)^{1/q}\\
&\le&2^{J(s+1/2-1/p)}(\sum_{j=0}^J||\beta_{j.}-\beta_{j.}^0||_2^q)^{1/q}\\
&\le&2^{J(s+1/2-1/p)}J^{\max(1/q-1/2,0)}||P_J\beta-P_J\beta^0||_2\\
&\le&n^{(s+1/2-1/p)/(2s+1)}\cdot c\tau_n
\end{eqnarray*}

Case 2: $1\le p<2$.

Since $||\beta_{j.}||_p\le 2^{j(1/p-1/2)}||\beta_{j.}||_2$ when
$1\le p<2$, the $B^s_{p,q}$ norm for $P_J\beta-P_J\beta^0$ can be
bounded as follows:
\begin{eqnarray*}
&&||P_J\beta-P_J\beta^0||_{B^s_{p,q}}\\
&\le&(\sum_{j=0}^J 2^{j(s+1/2-1/p)q}||\beta_{j.}-\beta_{j.}^0||_p^q)^{1/q}\\
&\le&(\sum_{j=0}^J 2^{jsq}||\beta_{j.}-\beta_{j.}^0||_2^q)^{1/q}\\
&\le&2^{Js}(\sum_{j=0}^J||\beta_{j.}-\beta_{j.}^0||_2^q)^{1/q}\\
&\le&2^{Js}J^{\max(1/q-1/2,0)}||P_J\beta-P_J\beta^0||_2\\
&\le&n^{s/(2s+2-2/p)}\cdot c\tau_n
\end{eqnarray*}

Summarizing the above two cases,
$||P_J\beta-P_J\beta^0||_{B^s_{p,q}}$ will be less than $\delta/2$
when $c$ is sufficiently small, and (\ref{temp}) is proved by
noticing $||P_J\beta||_{B^s_{p,q}}\le
||P_J\beta-P_J\beta^0||_{B^s_{p,q}}+||\beta^0||_{B^s_{p,q}}<B-\delta/2$

What is left is to lower bound
$\Pi_n(\sum_{j=0}^J\sum_k(\beta_{jk}-\beta_{jk}^0)^2\le c^2
\tau^2_n/\log n))$

Obviously the above prior probability is smallest when $\pi_j=1$
and thus the prior is a normal distribution. Let
$\delta_n^2=c^2\tau_n^2/\log n$ for simplicity of notation. If we
denote by $K=\sum_{j=0}^J j2^j\asymp (\log_2 n) n^{1/\alpha}$ the
total number of variables $\beta_{jk}, \mbox{ with } 0<j\le J$,
and let $\Delta=\exp\{-\sum_{j=0}^J\sum_k 2^{\alpha
j}(\beta_{jk}^0)^2\}$, $A=\{w_{jk}: 0\le j\le J, 0\le k\le 2^j-1,
||w||_2^2\le \delta_n^2\},$ we have
\begin{eqnarray*}
&&\Pi_n(\sum_{j=0}^J\sum_k(\beta_{jk}-\beta_{jk}^0)^2\le\delta_n^2)\\
&=&(\frac{1}{2\pi})^{K/2}\prod_{j=0}^J(2^{\alpha
j/2})^{2^j}\int_A\exp\{-\frac{1}{2}\sum_{j,k}2^{\alpha
j}(w_{jk}+\beta^0_{jk})^2\}\\
&\ge&\Delta(\frac{1}{2\pi})^{K/2}\prod_{j=0}^J(2^{\alpha
j/2})^{2^j}\int_A\exp\{-\sum_{j,k}2^{\alpha j}(w_{jk})^2\}\\
&\ge& \Delta(\frac{1}{2\pi})^{K/2}\prod_{j=0}^J(2^{\alpha
j/2})^{2^j}\frac{(\delta_n)^K\cdot\pi^{K/2}}{\Gamma(K/2)}\int_0^1u^{(K/2-1)}\exp\{-2^{\alpha
J}\delta_n^2\cdot u\}\,du \\
&=&\Delta(\frac{1}{2\pi})^{K/2}\prod_{j=0}^J(2^{\alpha
j/2})^{2^j}\frac{(\delta_n)^K\cdot\pi^{K/2}}{\Gamma(K/2)}(2^{-\alpha
J})^{K/2}(\delta_n^2)^{-K/2}\int_0^{2^{\alpha
J}\delta_n^2}u^{K/2-1}e^{-u}du\\
&\ge& 2^{-C\alpha JK}F(2^{\alpha J}\delta_n^2;K/2)\\
&\ge& e^{-Cn\epsilon_n^2} e^{K/2-2^{\alpha J}\delta_n^2}(2^{\alpha
J+1}\delta_n^2/K)^{K/2}(K/2)^{-1/2}\\
&\ge&e^{-Cn\epsilon_n^2}
\end{eqnarray*}
In the above we used Lemma 3 in \citet{shenwong01} as well as the
inequality $F(b;\alpha):=\frac{1}{\Gamma(\alpha)}\int_0^b
x^{\alpha-1}e^{-x}dx\succsim
e^{\alpha}e^{-b}b^{\alpha}\alpha^{-\alpha}\alpha^{-1/2}$ which
also appeared in that paper.  $\quad\Box$

If we use the sieve prior presented right before Theorem
\ref{maintheorem}, the same conclusion still holds.

\begin{theorem}\label{sievetheorem}
The result of Theorem \ref{maintheorem} is still true with the
sieve prior.
\end{theorem}
\textbf{Proof.} The entropy bound for condition (I) is unchanged.
With the same $J=\log_2 n/\alpha$, we have

\begin{eqnarray*}
&&\Pi_n^B(||\beta-\beta^0||_2^2\le\epsilon_n^2) \\
&\ge&\Pi_n(||\beta-\beta^0||_2^2\le\epsilon_n^2,
||\beta||_{B^s_{p,q}}<B\}\\
&\ge&\lambda_J\Pi_n^J(\sum_{j=0}^J\sum_k(\beta_{jk}-\beta_{jk}^0)^2\le\epsilon^2_n/2,
||P_J\beta||_{B^s_{p,q}}<B)\cdot\\
&&\Pi_n^J(\sum_{j=J+1}^\infty\sum_k(\beta_{jk}-\beta_{jk}^0)^2\le\epsilon^2_n/2)
\end{eqnarray*}

In the second probability above the event is actually
deterministic since $\beta_{jk}=0$ when $j>J$ under the prior
$\Pi_n^J$, and $\sum_{j>J,k}(\beta^0_{jk})^2\le\epsilon_n^2/2$. So
the probability of this term is $1$ and Lemma \ref{lemma:secterm}
is not needed.

For the first probability, the lower bound is exactly the same as
above. So the lower bound for the prior probability
$\Pi_n^B(||\beta-\beta^0||_2^2\le\epsilon_n^2)$ is bounded below
by $\lambda_J e^{-Cn\epsilon_n^2}$, and
$\lambda_J=n^{-\mu/\alpha}$ is obviously ignorable (can be
incorporated into the constant $C$ in $e^{-Cn\epsilon_n^2}$) in
this case. $\quad\Box$

If we focus on bounded functions only, then we can get the rate of
convergence for the posterior mean and posterior median.

\begin{corollary}
Consider the case where $||\beta^0||_2< 1$ and the prior on
$\beta$ is also renormalized to put mass $1$ on the set
$\{\beta:||\beta||_2<1\}$. Denote by $\hat{\beta}$ and
$\tilde{\beta}$ the posterior mean and posterior median
respectively. We have $||\hat{\beta}-\beta^0||_2=O(\epsilon_n)$
and $||\tilde{\beta}-\beta^0||_2=O(\epsilon_n)$ in probability.
\end{corollary}
\textbf{Proof.} First note that with slight modifications, Theorem
\ref{maintheorem} and Theorem \ref{sievetheorem} are still true
when the prior is constrained to unit $l_2$ balls.

The result for posterior mean is well-known
\citep{barron99,ghosal00} since the $l_2$ loss is bounded under
the current assumptions.

For posterior median, since the $l_2$ loss is now bounded and the
posterior probability $\Pi_n^B(||\beta-\beta^0||_2^2\ge
M\epsilon_n^2|X)$ converges to zero at least at the order
$\epsilon_n^2$ (implicit in the proof of \citet{ghosal00}, Theorem
2.1), we have $E||\beta-\beta^0||_2^2=O(\epsilon_n^2)$ in
probability, where the expectation is over the posterior
distribution of $\beta$. Then we use the simple fact that for any
random variable $X$, $E[X^2]\le a^2$ implies $|median(X)|\le 2a$.
This can be seen by $P(|X|>2a)\le E(X^2)/(4a^2)< 1/2$. Now
replacing $X$ by $\beta_{jk}-\beta^0_{jk}$, and summing over $j$
and $k$, we get the convergence rate for $\tilde{\beta}$.
$\quad\Box$

\section{Discussion}\label{sec:discussion}
Using the approach of \citet{ghosal00,ghosal07}, we have
investigated the convergence rate of the posterior distribution
for Gaussian white noise model in Besov spaces. Investigation of
posterior distribution rather than the Bayes estimators seems to
be more desirable from a philosophical and practical point of
view, since the posterior distribution can be directly utilized to
assess the uncertainty of the Bayesian inference. As shown in
\citet{abramovich04}, their Bayes factor estimator can achieve a
better rate of convergence (although it is still not optimal
within the whole range $1\le p<2$). Using the prior (\ref{prior})
we cannot hope to achieve this rate since it was shown in
\citet{abramovich04} that the posterior mean cannot achieve this
faster rate and the rate for the posterior distribution is no
faster than that of the posterior mean.

The loss function used in this investigation is the simplest $l_2$
loss. The extension to more general $l_p$ norm is left for further
research. The derived rate is the same as in \citet{abramovich04}
up to an extra $\log$ term and is suboptimal in the inhomogeneous
cases $1\le p<2$. Heavy-tailed distributions like double
exponential are successfully used in \citet{johnstone05} to
achieve better rates and it was argued that the implicit
thresholding in normal mixture are too heavy on high-resolution
levels. We believe optimal rates for posterior distribution are
achievable with similar heavy-tailed distributions.

% \appendix

% \section{}
% \label{}

% Bibliographic references with the natbib package:
% Parenthetical: \citep{Bai92} produces (Bailyn 1992).
% Textual: \citet{Bai95} produces Bailyn et al. (1995).
% An affix and part of a reference:
%   \citep[e.g.][Ch. 2]{Bar76}
%   produces (e.g. Barnes et al. 1976, Ch. 2).
\bibliographystyle{elsart-harv}
\bibliography{papers}

\end{document}